\documentclass[fleqn]{article}
\usepackage{amsfonts,amssymb}
\usepackage{latexsym, epsf}

\usepackage{times}

\newcommand{\R}{{\bf R}}

\newcommand{\forconf}[1]{}

\newcommand{\controls}{{\bf U}}

\newcommand{\states}{{\cal X}}

\newcommand{\rref}[1]{(\ref{#1})}

\newcommand{\kk}{{\cal K}}
\newcommand{\kl}{{\cal K}{\cal L}}
\newcommand{\ki}{{\cal K_\infty}}

\newcommand{\calr}{{\cal R}}

\newcommand{\norm}[1]{\left\Vert#1\right\Vert}

\newcommand{\abs}[1]{\left\vert #1 \right\vert}

\newcommand{\ve}{\varepsilon }

\newcommand{\E}{{\cal E}}

\newcommand{\B}{{\cal B}}

\newcommand{\A}{{\cal A}}

\def\edo{\end{document}}

\newtheorem{theorem}{Theorem}
\newtheorem{itlemma}{Lemma}[section]
\newtheorem{itproposition}[itlemma]{Proposition}
\newtheorem{itcorollary}[itlemma]{Corollary}
\newtheorem{itremark}[itlemma]{Remark}
\newtheorem{itdefinition}[itlemma]{Definition}
\newtheorem{itexample}[itlemma]{Example}
\newenvironment{lemma}{\begin{itlemma}\rm}{\end{itlemma}} 
\newenvironment{remark}{\begin{itremark}\rm}{\end{itremark}} 
\newenvironment{corollary}{\begin{itcorollary}\rm}{\end{itcorollary}}
\newenvironment{proposition}{\begin{itproposition}\rm}{\end{itproposition}}
\newenvironment{definition}{\begin{itdefinition}\rm}{\end{itdefinition}}
\newenvironment{example}{\begin{itexample}\rm}{\end{itexample}}
\def\bi{\begin{itemize}}
\def\ei{\end{itemize}}
\def\ben{\begin{enumerate}}
\def\een{\end{enumerate}}
\def \beq {\begin{eqnarray}}
\def \eeq {\end{eqnarray}}
\def \beqn {\begin{eqnarray*}}
\def \eeqn {\end{eqnarray*}}

\newcommand{\twoif}[4]{
\left\{ \begin{array}{ll}#1&#2\\#3&#4\end{array}\right.
}
\newcommand{\text}[1]{\hbox{\rm \ #1\ \/}}
\newcommand{\be}[1]{\begin{equation}\label{#1}}
\newcommand{\ee}{\end{equation}}
\newcommand{\bl}[1]{\begin{lemma}\label{#1}}
\newcommand{\br}[1]{\begin{remark}\label{#1}}
\newcommand{\bt}[1]{\begin{theorem}\label{#1}}
\newcommand{\bd}[1]{\begin{definition}\label{#1}}
\newcommand{\bp}[1]{\begin{proposition}\label{#1}}
\newcommand{\bc}[1]{\begin{corollary}\label{#1}}
\newcommand{\bfact}[1]{\begin{fact}\label{#1}}
\newcommand{\bex}[1]{\begin{example}\label{#1}}
\newcommand{\bem}[1]{\begin{example}\label{#1}}  
\newcommand{\ec}{\mybox\end{corollary}}
\newcommand{\efact}{\mybox\end{fact}}
\newcommand{\eex}{\mybox\end{example}}
\newcommand{\eem}{\mybox\end{example}}
\newcommand{\el}{\mybox\end{lemma}}
\newcommand{\ele}{\mybox\end{lemmaex}}
\newcommand{\er}{\mybox\end{remark}}
\newcommand{\et}{\qed\end{theorem}}
\newcommand{\ed}{\mybox\end{definition}}
\newcommand{\ep}{\mybox\end{proposition}}
\newcommand{\epr}{\end{proof}}
\newcommand{\bpr}{\begin{proof}}

\newcommand{\ecs}{\end{corollary}}
\newcommand{\eers}{\end{exercise}}
\newcommand{\eexs}{\end{example}}
\newcommand{\eems}{\end{example}}
\newcommand{\els}{\end{lemma}}
\newcommand{\eles}{\end{lemmaex}}
\newcommand{\ers}{\end{remark}}
\newcommand{\ets}{\end{theorem}}
\newcommand{\eds}{\end{definition}}
\newcommand{\eps}{\end{proposition}}
\newcommand{\qed}{\hfill \halmos} 
\newcommand{\mybox}{\hfill $\Box$} 
\newcommand{\comment}[1]{}
\newcommand{\obsolete}[1]{}
\newcommand{\halmos}{\rule{1ex}{1.4ex}}
\newenvironment{proof}{\noindent {\em Proof}.}
{\hspace*{\fill}$\halmos$\medskip}

\newcommand{{\rd}}{{\rm d}}
\newcommand{\rv}{{\rm v}}

\newcommand{\mix}{{\phi }}

\newcommand{\sqb}[1]{\left[ #1 \right]}
\newcommand{\cbrace}[1]{\left\{ #1 \right\}}

\newcommand{\ds}{\displaystyle}
\newcommand{\dist}{\mbox{dist}}

\newcommand{\uc}{{\bf u}}
\newcommand{\vc}{{\bf v}}

\newcommand{\forpaper}[1]{}

\newcommand{\maf}{{\cal X}}

\newcommand{\absa}[1]{{\abs{#1}}_{\A}}

\newcommand{\pna}{{\cal PN^A}}

\newcommand{\tn}[2]{{\bf TN_{#1}} {#2}}
\newcommand{\tube}{\Xi}
\renewcommand{\rv}{{\rm q}}
\newcommand{\pin}{\pi_N}
\newcommand{\pim}{\pi_M}
\renewcommand{\dist}{{\rm dist}}
\newcommand{\ab}[1]{{\langle #1 \rangle}}
\newcommand{\clos}{{\rm clos} \,}
\newcommand{\wc}{{\bf w}}

\begin{document}
\thispagestyle{empty} \setcounter{page}{0}

\title{Global Stabilization
   for Systems Evolving on Manifolds \footnote{
Running Head: Global Stabilization
   for Systems Evolving on Manifolds.
 Contact Author: Michael
Malisoff; malisoff@lsu.edu; Tel: (225) 578-1665;  Fax: (225)
578-4276.}}

\author{{}\\{}\\
 Michael
Malisoff\\ Department of Mathematics\\
Louisiana State University\\ Baton Rouge, LA 70803 USA\and{}\\{}\\Mikhail Krichman\\ ALPHATECH, Inc.\\
6 New England Executive Park\\ Burlington, MA 01803 USA
  \and {}\\Eduardo Sontag\\
  Department of Mathematics\\
Rutgers University\\  Hill Center-Busch Campus\\ New Brunswick, NJ
08903 USA}
\date{}

 \maketitle\newpage
 \begin{abstract}We show that any globally asymptotically controllable
 system on any smooth manifold can be globally stabilized by a
state feedback.  Since we allow discontinuous feedbacks, we
interpret the solutions of our systems   in the ``sample and
hold'' sense  introduced by Clarke-Ledyaev-Sontag-Subbotin (CLSS).
Our work generalizes the CLSS Theorem which is the special case of
our result for systems  on Euclidean space.  We apply our result
to the input-to-state stabilization  of systems on manifolds
relative to actuator errors, under small observation
noise.\medskip

\end{abstract}

\section{Introduction}
\label{sec1}

This note is devoted to the study of fully nonlinear systems
\begin{equation}
\label{systema} \dot x=f(x,u),\; \; x\in {\cal X},\; \; u\in
{\mathbf U}
\end{equation}
evolving on  arbitrary smooth manifolds ${\cal X}$ with inputs $u$
in general locally compact metric spaces  ${\mathbf U}$, where $f$
is locally Lipschitz in $x$ uniformly for $u$ in compact sets, and
jointly continuous in $(x,u)$.  We assume that (\ref{systema}) is
globally asymptotically controllable (GAC) to a given compact
weakly invariant nonempty set ${\cal A}\subseteq {\cal X}$; see
Section \ref{sec3} below for the definition of GAC for  systems on
manifolds.

It is natural to inquire about the relationship between the GAC
property for (\ref{systema}) and the existence of a feedback
$k(x)$ such that the closed-loop system
\begin{equation}
\label{systemb} \dot x=f(x,k(x)),\; \; x\in {\cal X}
\end{equation}
is globally asymptotically stable to ${\cal A}$. For the special
case where the system \rref{systema} evolves on ${\cal X}=\R^n$
and ${\cal A}=\{0\}$, this relationship has been well studied (see
\cite{B83,C92,MRS04, MS04, SS80}). For that case, it is now well
known that (\ref{systema}) does {\em not} in general admit a {\em
continuous} stabilizing $k(x)$ (see \cite{SS80, S79}). This
negative result can also be seen from Brockett's Criterion (see
\cite{B83, S99})
 which states that a
necessary condition for the existence of a continuous stabilizing
feedback $k(x)$ for \rref{systema} with $k(0)=0$ is that
$(x,u)\mapsto f(x,u)$ be open at zero;  see also \cite[pp.
252--255]{S98} for a simple direct proof of Brockett's result
using a homotopy. As a consequence, no totally  nonholonomic
mechanical system $\dot x=g_1(x)u_1+\ldots+g_m(x)u_m$ on $\R^n$
with $m<n$ and ${\rm rank}[g_1(0),\ldots, g_m(0)]=m$  is
stabilizable by a continuous state feedback (see \cite{S99}). On
the other hand, if \rref{systema} is GAC to $\A=\{0\}$ on $\R^n$,
then it can be
 stabilized by a
continuous {\em time varying } feedback $u=k(t,x)$ provided (i)
the system is completely controllable with no drift or  (ii) $n=1$
(see \cite{C92, SS80} and Remark \ref{rk3} below).

However, if we allow {\em dis\/}continuous feedbacks, then we have
the following positive result  from \cite{CLSS97}  known as the
Clarke-Ledyaev-Sontag-Subbotin (CLSS) Theorem: {\em If
(\ref{systema}) is GAC to ${\cal A}$ on $\maf=\R^n$, then there
exists a discontinuous feedback $k(x)$ for which (\ref{systemb})
is globally asymptotically stable to ${\cal A}$.}  Here and in the
sequel, `discontinuous' means `not necessarily continuous in the
state variable'. The discontinuous feedback $k(x)$ produces a {\em
discontinuous right-hand side} in \rref{systemb}, which requires a
more general interpretation of solutions that can be applied to
  discontinuous dynamics. In \cite{CLSS97}, this issue
is resolved by interpreting the trajectories of \rref{systemb} as
``sample and hold'' (a.k.a. CLSS) solutions (see Definition
\ref{pi-traj} below). The CLSS solution concept has been used
extensively in nonlinear control analysis and controller design
including the input-to-state stabilization of systems relative to
actuator errors under small observation noise (see \cite{MRS04,
MS04, S99} and Section \ref{sec5} below). For example, CLSS
solutions have been used to stabilize  nonholonomic systems such
as Brockett's Example which are not stabilizable by continuous
state feedbacks (see \cite{MRS04, MS04}).

On the other hand,  many important GAC systems evolve on manifolds
other than $\R^n$ (e.g., stabilization of rigid bodies on the Lie
group of rotations $SO_3$) and are therefore not tractable by the
CLSS Theorem.  In fact, if \rref{systema} is GAC to a singleton
$\A=\{p\}$ and admits a continuous stabilizing feedback $k(x)$,
then a theorem of Milnor (see \cite{M64}) implies that $\maf$ is
diffeomorphic to Euclidean space. This is because the existence of
$k(x)$ would imply the existence of a smooth control-Lyapunov
function on $\maf$ that could be viewed as a Morse function with a
unique (possibly degenerate) critical point, and manifolds
admitting such a Morse function are diffeomorphic to Euclidean
space (see \cite{S99}). Therefore, even if \rref{systema} is
holonomic, there may still be topological obstacles to continuous
global stabilization when $\maf\ne\R^n$.

Motivated by these considerations, this note will extend the CLSS
Theorem to GAC systems on general smooth manifolds ${\cal X}$,
proving the existence of a discontinuous feedback $k(x)$ rendering
(\ref{systemb}) globally stable to  ${\cal A}$ in the sense of
CLSS solutions. We follow the construction proposed in \cite{S99b}
which can be summarized as follows.  We first embed ${\cal X}$ as
a closed submanifold $g({\cal X})$ of a Eucldean space $\R^{k}$
for some $k$, e.g.,  using the Whitney Embedding Theorem. Then we
extend the system to all of $\R^{k}$ in
 such a way that (a) points outside $g({\cal X})$ can be
 controlled to a tubular neighborhood of ${\cal A}$ and (b) $g({\cal
 X})$ is invariant for the extended system.  We then apply the
 CLSS Theorem to the extended system on $\R^{k}$  to design
 our feedback $k(x)$.  The restriction of this feedback to
  $g(\maf)$ provides the desired stabilizer for the
 original system.

This note is organized as follows.  In Section \ref{sec2}, we
review CLSS solutions and the CLSS Theorem.  We introduce the
relevant definitions for stability on manifolds in Section
\ref{sec3}. In Section \ref{sec4}, we prove our Generalized CLSS
Theorem on the discontinuous stabilization of (\ref{systema}) on
smooth manifolds.  We illustrate our discontinuous feedback
constructions in Section \ref{illus}. We close in Section
\ref{sec5} by applying our results to  the input-to-state
stabilization  of GAC systems on Riemannian manifolds relative to
actuator errors under small observation noise.  This extends the
corresponding results \cite{MRS04, MS04, S89} on input-to-state
stabilization for systems evolving on Euclidean space.

\section{CLSS Theorem on Euclidean Space}
\label{sec2} In this section, we review the main definitions and
results from \cite{CLSS97} on the stabilization of GAC systems  on
Euclidean space.  Throughout this section, our state space is
${\cal X}=\R^n$. We  extend this material to systems on  smooth
manifolds in the next sections. We consider a system
(\ref{systema}) for which  $f$ is locally Lipschitz in $x$
uniformly for $u\in {\mathbf U}$ in compact sets, and jointly
continuous in $x$ and $u$.  Our input set ${\mathbf U}$ is a
locally compact metric space with a metric $d_{\mathbf U}$ and  a
distinguished element $0\in {\mathbf U}$, and we set
$|u|=d_{\mathbf U}(u,0)$ for each $u\in {\mathbf U}$. We let
${\mathbf {\cal U}}$ denote the set of all {\em controls} for
(\ref{systema}), i.e., the set of all
 measurable, locally essentially bounded functions $\uc : \R_{\ge 0} \to
 \controls$.  The essential supremum of any control $\uc\in {\mathbf {\cal U}}$
 is
 denoted by  $||\uc||$, and \[{\cal U}_N=\{\uc\in {\cal
U}: ||\uc||\le N\}\] for each $N>0$.
   Given $\xi\in \states$ and $\uc\in {\mathbf {\cal U}}$, the maximal
 trajectory of (\ref{systema}) for the control $\uc$ that satisfies
 $x(0)=\xi$ is denoted by $x(t,\xi,\uc)$ or simply by $x(t)$ when
 $\xi$ and $\uc$ are clear.  We say that $x(t)$
 is {\em well defined } provided it is defined for all $t\in
 \R_{\scriptscriptstyle \ge 0}:=[0,\infty)$.

Let $\A \subseteq \states$. We say that $\A$ is  {\em weakly
invariant (for \rref{systema})} provided there exists $N>0$ such
that for any $\xi \in \A$ there is a control $\uc\in {\cal U}_N$
such that the corresponding trajectory $x(t, \xi, \uc)$ is well
defined and stays in $\A$.  For example, $\A=\{0\}$ is weakly
invariant if $f(0,\bar a)=0$ for some $\bar a\in {\mathbf U}$.
More generally, $\A$ could be a periodic orbit we wish to
stabilize. We let $\abs p$ denote the Euclidean norm of any $p\in
\states$.  We let {\rm bd} (resp., {\rm clos}) denote the boundary
(resp.,  closure) operator, and we define the distance ${\rm
dist}({\cal N}, x)=\inf\{|p-x|: p\in {\cal N}\}$ for any subset
${\cal N}\subseteq\R^n$ and $x\in \R^n$. For any $x \in {\cal X}$,
we let $\absa{x}$ denote the distance from $x$ to $\A$.
 Therefore, $\absa{x} < \ve$ means  $x \in {\cal B}_\ve(\A):=\{p\in \R^n:
{\rm dist}(\A,p)< \ve\}$.

We next state two equivalent definitions of globally asymptotic
controllability.  First we state the well known definition from
\cite{MRS04, S99} in terms of comparison functions. We then
provide the original $\ve$-$\delta$ formulation which we
generalize to systems on manifolds in the next section. We use the
following
 comparison
function definitions from \cite{S99}.  A function $\alpha:
\R_{\scriptscriptstyle \geq 0} \to \R_{\scriptscriptstyle \geq 0}$
is said to be of class $\kk$ provided $\alpha$  is continuous,
strictly increasing, and satisfies $\alpha(0) = 0$; it is of class
$\ki$ provided it is also unbounded. We say $\alpha:
\R_{\scriptscriptstyle \geq 0} \to \R_{\scriptscriptstyle \geq 0}$
is of class
 $\cal N$ provided $\alpha$  is non-decreasing; and
of class   $\cal L$
  provided $\alpha(s)$ is decreasing to $0$ as $s \to +\infty$.
A function $\beta: \R_{\scriptscriptstyle \geq 0} \times
\R_{\scriptscriptstyle \geq 0} \to \R_{\scriptscriptstyle \geq 0}$
is said to be of class $\kl$ provided (a) $\beta(s, \cdot) \in
\cal L$ for every fixed $s$ and (b) $\beta(\cdot, t) \in \kk$ for
every fixed $t$. We write $\beta\in {\cal KL}$ to mean that
$\beta$ is of class ${\cal KL}$ and similarly for the other types
of comparison functions.

 \bd{kl_GAC} Let ${\cal A}\subseteq {\cal X}$ be compact, nonempty, and weakly
 invariant for \rref{systema}.
 We call \rref{systema}  {\em
globally asymptotically controllable (GAC) to $\A$ (on ${\cal
X}$)} provided there exist  $\beta \in \kl$ and $\sigma \in {\cal
N}$ such that for each $\xi \in \states$, there exists a control
$\uc$ with $\norm{\uc} \leq \sigma(\absa{\xi})$ such that $x(t,
\xi,\uc)$ is well defined and satisfies $ \absa{x(t,\xi,\uc)} \leq
\beta(\absa{\xi}, t)$ for all $t\ge 0$.\ed

 The following equivalent formulation of GAC has a natural
generalization  to systems  on manifolds;  see Definition
\ref{GAC} below. See  \cite{AS99} for  the equivalence of our GAC
definitions on $\R^n$.

\bd{tv_GAC_Rn} Let $\A \subseteq \states$ be  compact, nonempty,
and weakly invariant for \rref{systema}. We call \rref{systema}
{\em globally asymptotically controllable (GAC) to $\A$ (on ${\cal
X}$)} provided for all $\ve_1,\ve_2>0$ with $\ve_1<\ve_2$, we
have:

\begin{enumerate}\addtolength{\itemsep}{-0.25\baselineskip}
\item There exist $T=T(\ve_1, \ve_2) >0$ and  $\delta =
\delta(\ve_1)
>0$
 such that
for each $\xi \in \B_{\ve_2}(\A)$, there exists a control $\uc$
such that
\begin{itemize}
\item[(a)] $x(t, {\xi}, \uc)$ is well defined; \item[(b)] $x(t,
{\xi}, \uc) \in \B_{\ve_1}(\A)$ for all $t > T$; and \item[(c)] if
also ${\xi} \in \B_\delta(\A)$, then $\uc$ can be chosen so that
 $x(t,  {\xi}, \uc) \in \B_{\ve_1}(\A)$ for all $t \geq 0$.
\end{itemize}

\item For every positive number $\ve < \ve_2$,  there exists
$N=N(\ve)>0$ such that  if ${\xi}$ from 1. also satisfies
${\xi}\in {\cal B}_\ve({\cal A})$, then the control $\uc$ from 1.
can be chosen with $\uc\in {\cal U}_N$.\mybox
\end{enumerate}
\eds

\bd{obr_sv} A {\em feedback} for (\ref{systema}) is defined to be
any locally bounded function $k:  \states \to \controls$. \ed

In this note, we study the equivalence of  (open loop) asymptotic
controllability of  (\ref{systema}) and   the possibility of
stabilizing the system to a weakly invariant set $\A$ via a state
feedback.  The novelty of our work lies in its applicability to
systems on general smooth manifolds. Even for systems on $\R^n$,
it is often the case that a continuous stabilizing state feedback
does not exist (see \cite{MRS04, MS04, S99}).
 However, a discontinuous feedback is
always possible to construct, provided we use the
Clarke-Ledyaev-Sontag-Subbotin (CLSS) definition of a ``sample and
hold'' solution for a discontinuous dynamic.  We review  this
generalized solution notion next, following the notation from
\cite{MRS04, MS04}.

We define a {\em partition (of $\R_{\scriptscriptstyle \geq 0}$)}
to be any divergent sequence $ \pi: 0=t_0 < t_1 < t_2 < \ldots$
and we call
\[\ds\overline{\mathbf d}(\pi)=\sup_{i\ge 0}(t_{i+1}-t_i)\; \; \;
({\rm resp.,\ }\ds\underline{\mathbf d}(\pi)=\inf_{i\ge
0}(t_{i+1}-t_i))\] the {\em upper} (resp., {\em lower) diameter}
of the partition $\pi=\{t_o,t_1,t_2,\ldots\}$.

\bd{pi-traj} Let $k$ be a feedback for the system (\ref{systema}),
${\xi}\in {\cal X}$, and $\pi =\cbrace{t_i}_{i \geq 0}$ be a
partition. The {\em $\pi$-trajectory} \[t\mapsto
x_\pi(t,{\xi},k)\] for (\ref{systema}), ${\xi}$, $\pi$, and $k$ is
defined to be the continuous function obtained by recursively
solving \[\dot x(t)=f(x(t),k(x(t_i)))\] from the initial time
$t_i$ up to the maximal time \be{si}s_i=\max\{t_i,\sup\{s\in
[t_i,t_{i+1}]: x(\cdot) {\rm\ is\ defined\ on\ } [t_i,s)\}\},\ee
where $x(0)={\xi}$.\footnote{The continuity requirement for
$x(\cdot)$ amounts to stipulating that the final value on the
previous subinterval is used as the initial value at the next
subinterval.} The domain of $x_\pi(\cdot,{\xi},k)$ is  $[0, t_{\rm
max})$, where
\[t_{\rm max}=\inf\{s_i: s_i<t_{i+1}\}.\]  We call
$x_\pi(\cdot,{\xi},k)$ {\em well defined } provided $t_{\rm
max}=+\infty$. \ed

The $t_i$  argument in the maximum \rref{si}
 is needed to allow the
possibility that $x(\cdot)$ is not defined at all on $[t_i,
t_{i+1}]$ in which case the supremum in \rref{si} alone would by
definition give $-\infty$. The following notion of (global)
stabilization for (\ref{systema}) was introduced in~\cite{CLSS97}:
\bd{s-stab} A feedback $k: \states \to \controls$ is said to {\em
s-stabilize} the system~\rref{systema} to ${\cal A}$ provided for
each pair $(r,R)$ with  $0 < r < R$ there exist $M = M(R) > 0$
with $\lim_{R \to 0} M(R) = 0$,
 $\delta = \delta(r, R) > 0$, and $T = T(r, R) > 0$ such that,
for every  $\pi$ with $\overline{\bf d}(\pi)< \delta$ and
${\xi}\in {\cal B}_R({\cal A})$, the $\pi$-trajectory $x(\cdot)$
for  \rref{systema}, initial value $\xi$,  partition $\pi$, and
feedback $k$ is well defined and satisfies (a) $\absa{x(t)} \leq
r$ for all $t\ge T$ and (b) $\absa{x(t)} \leq M(R)$ for all $t\ge
0$.\mybox \eds The following {\bf result to be generalized} was
shown in \cite{CLSS97} for ${\cal A}=\{0\}$ but can be shown for
 our general compact, nonempty, weakly invariant set ${\cal
A}\subseteq {\cal X}=\R^n$ by  similar arguments (e.g., using the
existence results from \cite{KT04} for locally Lipschitz Lyapunov
functions for GAC systems and any compact set ${\cal
A}$):\bt{babyCLSS} If  (\ref{systema}) is GAC to ${\cal A}$ on
${\cal X}=\R^n$, then it admits a feedback that s-stabilizes the
system to ${\cal A}$. \et

The preceding result is called the {\em CLSS Theorem}.  Our main
contribution is a Generalized CLSS Theorem for systems on smooth
manifolds and is the subject of the next two sections.  We provide
related results on input-to-state stabilization  on Riemannian
manifolds in Section \ref{sec5}.

\section{Stabilization on Manifolds}
\label{sec3} We again consider the  system~\rref{systema} but we
assume from now on that the state space ${\cal X}$ for the system
is an arbitrary smooth (i.e., $C^\infty$)
(second countable) manifold. 
Controls $\uc$, as before, are measurable, locally essentially
bounded functions $\uc: \R_{\ge 0}\to\controls$.  We assume
\be{ti_f_on_m} f : \maf\times\controls \to T_x (\maf): (x, u)
\mapsto f(x, u) \ee is locally Lipschitz in $x$ and jointly
continuous in $x$ and $u$; that is
\[f(x, u) = \sum_i a_i(x, u)\frac{\partial}{\partial x_i}\]
where each $a_i:\maf\times {\mathbf U}\to \R$ is locally Lipschitz
in $x$ uniformly for $u$ in compact sets and jointly continuous,
and $T_x(\maf)$ is the tangent space to $\maf$ at $x$.
 We define the solutions $x(t,\xi,\uc)$ of (\ref{systema}) as
before.
 We next generalize Definition  \ref{tv_GAC_Rn} for GAC to
 manifolds.

Let $\A$ be a compact, nonempty, weakly invariant subset of $\maf$
for \rref{systema}, and let
 $\pna$  denote the set of all open precompact subsets of
 $\maf$
containing $\A$. To extend the GAC definition to manifolds,
 we  simply replace the $\ve$-neighborhoods of
$\A$ from Definition \ref{tv_GAC_Rn} with arbitrary sets in $\pna$
as follows:
 \bd{GAC}  We say that \rref{systema}
 is {\em globally
asymptotically controllable (GAC) to $\A$ (on $\maf$)} provided:
\begin{enumerate}
\item Given any  $\E_1, \E_2 \in \pna$ with $\E_1 \subseteq \E_2$,
 there exist $T=T(\E_1, \E_2) >0$ and     $\Delta = \Delta(\E_1) \in
\pna$
 such that
for every $\xi \in \E_2$ there exists a control $\uc$ such that
\begin{itemize}
\item[(a)] $x(t, \xi, \uc)$ is well defined; \item[(b)] $x(t, \xi,
\uc) \in \E_1$ for all $t > T$; and \item[(c)] if also $\xi \in
\Delta$, then $\uc$  can be chosen so that
 $x(t, \xi, \uc) \in \E_1$ for all $t \geq 0$.
\end{itemize}
\item For every set ${\cal N} \in \pna$,  there exists $N=N({\cal
N})>0$ such that if $\xi$ from 1. also satisfies $\xi \in \cal N$,
then the control $\uc$ from 1. can be chosen with $\uc\in {\cal
U}_N$. \mybox
\end{enumerate}
\eds

We assume throughout this section that our dynamic $f$ is GAC to
${\cal A}$. Since our definitions of feedback and $\pi$-trajectory
from Section \ref{sec2} do not depend on the structure of the
state space ${\cal X}$, they remain valid for systems  on
manifolds. We extend the definition of an s-stabilizing feedback
to manifolds as follows:

\bd{s-stab-maf} A feedback $k : \states \to \controls$ is said to
{\em s-stabilize} the system \rref{systema} to ${\cal A}$ provided
the following hold for all sets $\calr_1$, $\calr_2 \in \pna$ with
$\calr_1 \subseteq \calr_2$:
\begin{enumerate}
\item There exist a set ${\cal M} = {\cal M}(\calr_2)\subseteq
{\cal X}$ and
  numbers $\delta = \delta(\calr_1, \calr_2) >0$ and $T = T(\calr_1, \calr_2)$ such that,
 for any partition $\pi$ with $\overline{\bf d}(\pi) < \delta$
 and  any ${\xi}$ in $\calr_2$,
the $\pi$-trajectory $x(\cdot)$ for \rref{systema}, the initial
state $\xi$, and the  feedback $k$ is well defined and satisfies:
(a) $x(t) \in \calr_1$ for all $t \geq T$ and (b) $x(t) \in \cal
M$ for all $t \geq 0$.  \item[2.] For each set $\cal E\in\pna$
there exists $\cal D\in \pna$ such that if $\calr_2 \subseteq \cal
D$, then the set $\cal M$ in 1. can be chosen so that $\cal M
\subseteq E$.\mybox \end{enumerate}\eds

Our goal is to show that the CLSS Theorem remains true on any
smooth manifold $\maf$. To this end, we follow the strategy
outlined in \cite{S99b} which can be summarized as follows. We
first embed the   state space  manifold $\maf$ into some Euclidean
space  $\R^{k}$ (e.g., using the Whitney Embedding Theorem). Then
we extend the dynamic
 to  all of  $\R^{k}$ in such a way that (a) the
system is asymptotically controllable to a tubular neighborhood of
$\A$ and (b) $\maf$ is a strongly  invariant set under the
extended system (see Lemma \ref{invariance}). Next we apply the
CLSS Theorem to the extended system. Thus, we get an s-stabilizing
feedback on $\R^{k}$. When restricted to $\maf$, this feedback
will s-stabilize ~\rref{systema} to ${\cal A}$.

To make this construction precise, we use the following
definitions and facts from differential topology (see \cite{Bo03,
Br93}). The following is known as the Whitney Embedding Theorem
(see \cite[p.92]{Br93}):

\bl{WET} If ${\cal X}$ is an $n$-dimensional smooth manifold, then
there exists an
 embedding $g: {\cal X} \to \R^{2n+1}$ for which $g({\cal X})$ is a
 submanifold, and a
 closed subset, of
 $\R^{2n+1}$.
\el

By Lemma \ref{WET}, we can assume that our state space
 $\states$ is a smooth submanifold of $\R^k$  with $\maf\subseteq\R^k$ closed. The {\em
normal bundle} $\tube({\cal X})$ of ${\cal X}$ in $\R^k$ is
defined by \[\tube({\cal X})\; \; =\; \; \cbrace{\langle x, \rv
\rangle \in {\cal X} \times \R^k\,:\, \rv \bot T_x({\cal X})}.\]
We define the projections $\pim : \tube({\cal X}) \to {\cal X}$ by
$\pim (\langle x, \rv \rangle) = x$ and $\pin: \tube({\cal X}) \to
\R^{k}$ by $\pin (\langle x, \rv \rangle) = \rv$, and
$\theta(\langle x,q\rangle):=x+q$. For each smooth function
$\omega: {\cal X} \to \R_{>0}$, the {\em $\omega$-tube}
$\tube(\maf, \omega)$ is defined by
\[\tube(\maf, \omega)=\cbrace{\langle x, \rv \rangle \in \tube(\maf)
\,:\, \abs \rv < \omega(x)}.\] The next result is known as the
Tubular Neighborhood Theorem. \bl{TNT} Let $\maf$ be a closed
submanifold of $\R^k$. There exists a smooth function
$\omega:{\cal X}\to \R_{>0}$ such that $\theta:\tube(\maf, \omega)
\to \R^k: \langle x, \rv \rangle \mapsto x + \rv$ is a
diffeomorphism onto an open neighborhood of $\maf$ in $\R^k$. \el

In particular, $\tube(\maf, \omega)$ is an open subset of
$\R^k\times \R^k$.
   Pick  functions $\omega$ and $\theta$ as in Lemma~\ref{TNT} for
our state space manifold $\maf$.
 Since $\A\subseteq\maf$ is
compact and $\omega$ is continuous, $\omega$ attains its minimum
on $\A$. Let
\[\ve= \frac{1}{2}\min_{x \in \A} \omega(x)
\]
and for each set
 ${\cal A}_1\subseteq \maf$, define
\[\renewcommand{\arraystretch}{1.3}\begin{array}{l}
\tube({\cal A}_1, \ve)\; \; = \; \;  \cbrace{\langle x, \rv
\rangle \in \tube(\maf) \,:\, x \in {\cal A}_1, \abs \rv < \ve},\\
\tube({\cal A}_1, \omega)\;  =\;   \cbrace{\langle x, \rv \rangle
\in
\tube(\maf) \,:\, x \in {\cal A}_1, \abs \rv < \omega(x)},\\
\tn{\ve}{{\cal A}_1} \; \; =\; \;  \theta (\tube({\cal A}_1,
\ve)), \; \; \; \tn{\omega}{{\cal A}_1} \; \; =\; \;  \theta
(\tube({\cal A}_1, \omega)).\end{array}\]

Notice that if $\ve\le \omega(x)$ for all $x\in \A_1$,  then
$\tn{\ve}\A_1\subseteq \tn{\omega}\A_1\subseteq \tn{\omega}\maf$.
Also, ${\rm clos}\tn{\ve}\A$ is a compact subset of
$\tn{\omega}\maf$.
  Next consider the system \be{ext1} \dot x =
f(x, u),\; \; \; \dot \rv = \rv \, v,\; \; \; \langle
x,q\rangle\in \maf\times \R^k,\; \; \langle u,v\rangle \in
{\mathbf U}\times \R, \ee whose (maximal) solution for the
controls $\langle \uc, \vc\rangle$ starting at $\langle\xi,
q_0\rangle$ we denote by $\langle x(t, \xi, \uc), \rv(t, q_0, \vc)
\rangle$,
 or  by $\langle x(t), \rv(t) \rangle$ for brevity.
If, for some initial state $\ab{x(0), \rv(0)}$ and controls
$\langle \uc, \vc\rangle$, the  trajectory $\ab{x(t), \rv(t)}$ of
\rref{ext1} stays  in $\tube({\cal X}, \omega)$, then $y(t, y_0,
\uc, \vc) = \theta (\ab{x(t), \rv(t)})$
 is the corresponding trajectory of
 \be{ext2} \dot y = f_1(y, u, v) := f(\pim(\theta^{-1}(y)),
u) + \pin(\theta^{-1}(y))\, v,\; \; y\in \tn \omega \maf,\; \;
\langle u, v\rangle\in {\mathbf U}\times \R
 \ee  with the initial value $y_0 = y(0) = \theta (\ab{x(0),
\rv(0)})$. We denote this solution by $y(t)$ when no confusion
would result. We also omit  the $\theta^{-1}$ inside the
projections $\pi_N$ and $\pi_M$ in the sequel to simplify our
notation.  We (discontinuously) extend $f_1$ to $\R^k$ by defining
it to be zero outside $\tn \omega \maf$.

Next, we extend our GAC system \rref{systema} to  all of $\R^k$ as
follows. Let ${\cal X}^\sharp\subseteq\R^k$ be any  closed set
contained in $\tn{\omega}\maf$ and containing $\maf$ in its
interior. Let $C_\omega\subseteq\R^k$ be any open set such that
the following holds: \[\R^k\setminus \tn{\omega}\maf \subseteq
C_\omega\subseteq {\rm clos}\, C_\omega\subseteq
\R^k\setminus\maf^\sharp.\]  Then ${\rm bd} C_\omega
\subseteq\tn{\omega}\maf$.
 Let
$\mix:\R^k \to [0, 1]$ be any smooth function such that
\be{mix}\mix(z) = \twoif{1} { z \in \maf^\sharp}
                  {0} { z \in {\rm clos}\,
C_\omega}\ee which exists by a well known separation result (e.g.,
\cite[Exercise V.4.5]{Bo03}). Now define a system \be{ext3} \dot z
= f_2(z, u, v, w)
    := f_1(z, u, v) \mix(z) + (1-\mix(z))\, w,\; \; z\in \R^k,\; \;
\langle u, v, w\rangle  \in \controls \times \R \times \R^k, \ee
whose (maximal) solution starting at $z_o$ for  given controls
$\langle \uc, \vc, \wc\rangle$ we denote by $z(t, z_0, \uc, \vc,
\wc)$. Since $\phi\equiv 0$ in $C_\omega$, we know $f_2$ is
locally Lipschitz in $z\in \R^k$.  We use the following elementary
observation:

\bl{invariance} Any trajectory $z(t)$ for $f_2$ starting at a
point $\eta\in \maf$ remains in $\maf$ on its domain of definition
and therefore is a trajectory of $f$.  In other words, $\maf$ is
strongly invariant for $f_2$.\el
\begin{proof}
Since  $\langle x,0\rangle \in \Xi({\cal X},\omega)$ for all $x\in
\maf$, the
 uniqueness property for solutions of
\rref{ext2} in $\tn{\omega}{\maf}$ implies that all trajectories
of $f_1$ starting in ${\cal X}$ remain in ${\cal X}$ and so are
trajectories of $f$. On the other hand, trajectories $z(t)$ of
$f_2$ starting in ${\cal X}$ are also trajectories of $f_1$ while
they are in $\maf^\sharp$ (by our choice \rref{mix} of $\phi$),
since $f_1$ and $f_2$ agree on $\maf^\sharp$. By the uniqueness
property for trajectories of $f_1$, $z(t)$ therefore cannot enter
$\maf^\sharp\setminus\maf\subseteq \tn{\omega}\maf$ and so stays
in ${\cal X}$. Hence $z(t)$ is a trajectory of $f_1$, and also for
$f$.\end{proof}

\noindent The preceding lemma forms the basis for our Generalized
CLSS Theorem in the next section.

\section{CLSS Theorem on Manifolds}
\label{sec4} In this section, we prove  the  following {\em
Generalized  CLSS Theorem} for any  smooth manifold  ${\cal X}$
and any compact, nonempty, weakly invariant set ${\cal A}\subseteq
{\cal X}$ for \rref{systema}:

\bt{CLSS} If  (\ref{systema}) is GAC to ${\cal A}$ on the manifold
${\cal X}$, then it admits a feedback that s-stabilizes the system
to ${\cal A}$.\et

\noindent This will follow from the following key lemma:
\bl{ext3_is_GAC} If the system \rref{systema} is GAC to $\A$ on
$\maf$, then the system~\rref{ext3} is GAC to $\clos \tn{\ve}{\A}$
on $\R^k$. \el

We begin by proving Lemma \ref{ext3_is_GAC}.  Fix $z_0 \in \R^k$,
a precompact open set ${\cal B}$ containing ${\rm clos}
\tn{\ve}{\cal A}$, and an open set ${\cal A}_1\in {\cal PN}^{\cal
A}$ such that $\tn{\ve}{\cal A}_1\subseteq {\cal B}$.
 Assume first that $\eta:=z_0\in \tn{\omega}{\maf}$. Since~\rref{systema} is GAC to $\A$, we
can find a control $\uc: [0, \infty) \to \controls$ and constants
$T_1>0$
 and $p_1 > 0$ with
 $||\uc|| < p_1$ such that the trajectory
$x(t)=x( t, \pim(\eta), \uc)$ of \rref{systema}  is well defined
and satisfies $x(t)\in {\cal A}_1$ for all $t\ge T_1$.   This
gives a compact set $\bar B \subseteq \maf$ containing ${\cal A}$
such that $x(t) = x(t, \pim(\eta), \uc) \in \bar B$ for all $t\ge
0$.

Since $\omega$ is positive and  smooth on $\maf$, there exist
positive values
\begin{eqnarray}\label{bd} p_2= 1+\max_{x \in \bar B} \abs{\nabla
\omega(x)},\; \; \; p_4 = \min_{x \in\bar
B}\omega(x)\end{eqnarray} and   $p_3 > 0$ such that
 $\abs {f(x, u)} < p_3$ for all $x \in \bar B$ and $\abs u < p_1$.
Then $p_4/2\le \ve$, and

\be{efbound} \abs{\frac{d}{dt}\omega(x(t))} = \abs{\nabla
\omega(x(t)) \cdot f( x(t), \uc(t))}  \leq p_2 p_3\; \; \; {\rm
for\ \ almost\ \ all\ \ } t\ge 0.\ee In other words, $p_2 p_3$ is
an upper bound on the rate of change of the width $\omega(x(t))$
of $\tn{\omega}{\maf}$, as we move along the trajectory $x(t)$.
Hence, to ensure that our stabilizing trajectory of \rref{ext3}
starting in $\tn{\omega}{\maf}$ stays there, we must design a
control $\vc$  so that the solution  of \rref{ext2} is  pushed
towards $\maf$ faster than $p_2 p_3$.

Since we assumed $\eta \in \tn{\omega}{\maf}$, we have
$\langle\pi_M(\eta),\pi_N(\eta)\rangle\in \tube({\cal X},\omega)$
and therefore
 $\abs{\pin(\eta)} < \omega(\pim(\eta))$.  Define $\vc : [0, +\infty) \to \R$
 by

\be{vchoice}\vc(t) =\twoif { - \frac{p_2 p_3}{p_4/4}} {t \in
\sqb{0, T_2}} {0}{t > T_2},\; \; \; {\rm where}\; \; \; T_2 =
\max\cbrace {0, {\frac{\abs{\pin(\eta)} - p_4/4}{p_2 p_3}}}. \ee
Let $q(t)$ be the solution of $\dot q=\vc q$ starting at
$\pi_N(\eta)$.  Set $y(t)=x(t)+q(t)$, where $x(t)=x(t,\pi_M(\eta),
\uc)$ is defined above;  then $y(t)$ has domain $[0,\infty)$, and
$\langle x(t),q(t)\rangle$ is a solution of \rref{ext1} on
$[0,+\infty)$. We next define $t'=\inf\{t\ge 0: \langle
x(t),q(t)\rangle \in {\rm bd}\tube({\cal X},\omega)\}$, so
$\langle x(t),q(t)\rangle\in \tube(\maf, \omega)$ on $[0,t')$.
 We show
that $t'=+\infty$.  This will show that $y(t)$ is a solution of
\rref{ext2} on all of $\R_{\scriptscriptstyle \geq 0}$.  To this
end, first note that:
\begin{itemize}
\item[(i)] Since the direction of $\vc(t) q(t)$ is always opposite
to that of $q(t)$ whenever $\vc(t) \neq 0$, the function
$\abs{q(t)}$ is non-increasing on $\R_{\ge 0}$.

 \item[(ii)] At all
points  $t\in [0,T_2]$ for which $\dot x(t)$ exists and
$\abs{q(t)} \ge p_4/4$, the following holds:
\begin{equation}
\label{eq1:ext3_GAC} \frac{d}{dt}|q(t)|\; =\;
-\frac{p_2p_3}{p_4/4}|q(t)|\; \le \;
 -p_2 p_3 \;  \le \;
 -\abs{\frac{d}{dt}\omega(x(t))}\; \le \; \frac{d}{dt}\omega(x(t)).
 \end{equation}\end{itemize}
By separately considering the case where $|q(t)|$ stays above
$p_4/4$ on $[0,T_2]$ and using \rref{vchoice}-\rref{eq1:ext3_GAC},
one can easily check that $|q(T_2)| \le p_4/4$; this inequality is
clear if $|q(t)|$ ever goes below $p_4/4$ on $[0,T_2]$, by (i).
Hence, $|q(t)|\le p_4/4$ for all $t\ge T_2$, by the choice of
$\vc$.
 Similarly, we
can use \rref{eq1:ext3_GAC}, the definition  of $p_4$, and the
fact that $\abs{\pin(\eta)} < \omega(\pim(\eta))$ to verify that
\be{verify} \abs{q(t)} < \omega(x(t))\; \; \; \forall t\ge 0.\ee
Suppose that $t'<\infty$.  Then $\langle x(t'), q(t')\rangle\in
{\rm bd}\tube(\maf, \omega)$.
Since $\tube({\cal X})$ is closed and $\langle x(t),q(t)\rangle
\in \tube({\cal  X})$ on $[0,t')$, it follows from \rref{verify}
that $\langle x(t'),q(t')\rangle \in \tube({\cal  X},\omega)$,
contradicting the openness of $\tube({\cal  X},\omega)$.
 It follows that $t'=+\infty$, so
the solution $y(t):=y(t, \eta, \uc, \vc)$ of the
system~\rref{ext2} maps all of $\R_{\ge 0}$ into
$\tn{\omega}{\maf}$.

Finally, we define a control $\wc: [0, \infty) \to \R^k$ by
\be{wchoice} \wc(t)= f_1( y(t, \eta, \uc, \vc), \uc(t), \vc(t)).
\ee The control $\mathbf w$ cancels the effect of $\mix$ in
\rref{ext3} for states in
$\tn{\omega}{\maf}$.  In fact,
\[f_2(y(t, \eta, \uc, \vc), \uc(t), \vc(t), \wc(t))\equiv f_1(y(t,
\eta, \uc, \vc), \uc(t), \vc(t)),\] hence $y(t, \eta, \uc, \vc)
\equiv z(t,\eta, \uc, \vc, \wc)$.  By our choices  of $T_1$,
$p_4$, and $\vc$, we have (a) $x(t)=\pi_M(y(t,\eta,\uc,\vc))\in
{\cal A}_1$ for all $t\ge T_1$ and (b) $\abs{\pin( y(t,\eta, \uc,
\vc))}< p_4/2 \le \ve$ for all $t> T_2$.  It therefore follows
that  \[z(t, \eta, \uc, \vc, \wc) =
  y(t, \eta, \uc, \vc) \in  \tn{\ve}{\A_1} \subseteq \B\; \;
 \forall t >  T,\] where $T:= \max\cbrace{T_1,
T_2}$.  This shows the asymptotic controllability of
 (\ref{ext3}) to our arbitrary neighborhood ${\cal B}$ of  ${\rm clos}\tn{\ve}{\cal A}$ from any initial
 value in $\tn{\omega}{\cal X}$.  We next show that this
 controllability property holds from initial values outside $\tn{\omega}{\cal
 X}$ as well.

Assume then that $z_0 \not\in \tn{\omega}{\maf}$, so $z_0\in
C_\omega$. We reduce to the case where  the initial value is in
$\tn{\omega}{\maf}$. Let $p_5 = \dist(\maf, z_0)$ and let $\eta_1
\in \maf$ be  such that $\abs{\eta_1 - z_0} = p_5$. Define $\bar
\wc$ and $z(t)$ by \begin{equation} \label{wchoice2}\label{zform}
 \bar \wc(t) =
\frac {\eta_1 - z_0}{\abs{\eta_1 - z_0}} \;\;\forall \, t \ge 0;\;
\; \; \; \;  z(t)=z_o+t\frac{\eta_1-z_0}{|\eta_1-z_0|},\; \; \;
0\le t\le \hat t:=\inf\{t\ge 0: z(t)\in {\rm bd}\, C_\omega\}.
\end{equation}Then $z(t)$ is a solution of \rref{ext3} starting
at $z_0$ for any controls $\langle \uc,\vc\rangle$ and the choice
$\wc=\bar \wc$, and $z(t)\in C_\omega$ on $[0,\hat t)$.
   Also, $0<\hat t\le |\eta_1-z_0|$, since if it were the case
   that $\hat t> |\eta_1-z_0|$, then
setting $t=|\eta_1-z_0|$ in \rref{zform} would give
$z(t)=\eta_1\in {\cal X}\cap C_\omega$.  This would contradict the
fact that $C_\omega\subseteq \R^k\setminus\maf$. We conclude in
particular that $\hat t<\infty$, so  $\eta := z(\hat t, z_0, \uc,
\vc, \bar \wc) \in {\rm bd}C_\omega\subseteq \tn{\omega}{\maf}$.
 For our
precompact open set ${\cal B}$, we now construct the controls $
\uc$ from the controllability of~\rref{systema}, $\vc$ as
in~\rref{vchoice}, and $\wc$ as in~\rref{wchoice}, driving this
choice of $\eta$ to ${\cal B}$. Let $\uc^\sharp$ and $\vc^\sharp$
be the concatenations of the zero functions on $[0,\hat t)$,
followed by $\uc$ and $\vc$, respectively.  Let $\wc^\sharp$ be
the concatenation of  $\bar \wc$ on $t \in [0, \hat t)$ from
\rref{wchoice2},
 followed by
$\wc$ from \rref{wchoice} for $t \geq \hat t$.  The control vector
$\langle \uc,\vc,\wc\rangle$ for \rref{ext3} drives $\eta$ to $\B$
in time $T$, so $z(t, z_0, \uc^\sharp, \vc^\sharp, \wc^\sharp) \in
\B$ for all $t
> \hat t + T$. Since $T$ and $\hat t$
are locally bounded functions of $z_o$ and ${\cal B}$, we conclude
that Conditions 1 (a)-(b) from  the GAC definition  hold for
\rref{ext3} and the attractor ${\rm clos}\tn{\ve}{{\cal A}}$.

To establish Condition 1 (c) of the GAC definition for
\rref{ext3}, fix any precompact open set $E\subseteq \R^k$
containing $\clos \tn{\ve}{\A}$. We can find an open  set $\E_1
\in \pna$ and $\ve'>\ve$ such that $\ve'<\omega(x)$ for all $x\in
\E_1$, and such that $\tn{\ve'}{\E_1} \subseteq E$. Next we find a
set $\Delta\in\pna$ as in Condition 1 (c) of Definition~\ref{GAC}
for the GAC system~\rref{systema}, corresponding to ${\cal E}_1$.
It follows that ${\rm clos} \tn{\ve}{{\cal A}}\subseteq
\tn{\ve'}{\Delta}$. By reducing
 $\Delta$, we can assume $\ve'<\omega(x)$ for all $x\in
 \Delta$, and therefore
 $\tn{\ve'}\Delta\subseteq  \tn{\omega}\Delta\subseteq
 \tn{\omega}\maf$. We show
that if $z_0 \in D:= \tn{\ve'}{\Delta}$, then $z_0$ can be driven
to ${\cal B}$ using the  system \rref{ext2} and the vector of
controls $\langle \uc, \vc\rangle$ as defined above, hence also by
the extended system \rref{ext3}, while being kept inside
$\tn{\ve'}\E_1\subseteq E$ for all $t\ge 0$.

 Let $z_0 \in D$.   Since $\pi_M(z_o)\in \Delta$, we can
 arrange (by the choice of $\Delta$) that
 $\uc$ is such that
 $x(t, \pim(z_0), \uc) \in \E_1$ for all $t \geq 0$.
Next, we construct  $\vc$ defined by \rref{vchoice} for the
initial state   $z_0\in \tn{\ve'}\Delta\subseteq \tn{\omega}{\cal
X}$. By (i), we know $t\mapsto |\pin(y(t, z_0, \uc, \vc))|$ is
non-increasing. Thus,
 for all $t \ge 0$, we get
$|\pin(y(t, z_0, \uc, \vc))| \le |\pin(z_0)| < \ve'$ and therefore
also
 $y(t, z_0, \uc, \vc) \in \tn{\ve'}{\E_1} \subseteq E$,
proving Condition 1 (c) from the GAC definition for
system~\rref{ext3}.

It remains to check that the concatenated controls $\uc^\sharp$,
$\vc^\sharp$, and $\wc^\sharp$ we  constructed above satisfy the
boundedness requirement from Condition 2 of the GAC definition.
That is, we need to check that $||\langle \uc^\sharp, \vc^\sharp,
\wc^\sharp\rangle ||$ is a locally bounded function of the initial
state $z_0$. To do this, first note that the boundedness
requirement on $\uc^\sharp$ is satisfied because \rref{systema} is
assumed to be GAC to $\A$ on $\maf$. Next, $||\vc^\sharp|| \leq
p_2 p_3/(p_4/4)$, and (letting $\uc$ be the second part of the
concatenation $\uc^\sharp$ and similarly for $\vc$, as before)
\begin{eqnarray}
\norm{\wc^\sharp} \leq 1+
 {\rm ess}\sup_{t\ge 0}\left\{ |f(\pim(y(t, \eta,  \uc,  \vc)),
 \uc(t))|+|\pi_N(y(t,\eta,  \uc, \vc))|
 \frac{p_2 p_3}{p_4/4}\right\}.\nonumber \end{eqnarray} Here
$f(\pim(y(t, \eta,  \uc,  \vc)),
 \uc(t))$ stays bounded  because
 (a) $\pi_M(y(t,\eta,\uc, \vc)) =x(t,\pi_M(\eta),
 \uc)
 \in \bar B$ for all $t\ge
 0$ and  (b) $\bar B$ and the $p_i$'s are  locally bounded functions of the
 state $\eta=z(\hat t,z_o, \uc^\sharp,\vc^\sharp,\wc^\sharp)$.
 Also,
  $|\pi_N(y(t,\eta,\uc,\vc))|$ stays bounded because it
 decreases from $|\pi_N(\eta)|$.
  Hence, Condition 2 of the GAC definition holds.  This
 completes the proof of Lemma \ref{ext3_is_GAC}.

Finally, we prove  Theorem \ref{CLSS}.  The preceding argument
applied to $\eta\in {\rm clos} \tn{\ve}{\cal A}$ (with $\uc$
chosen so that $x(t,\pi_M(\eta), \uc)\in {\cal A}$ for all $t\ge
0$, which exists by the weak invariance of ${\cal A}$) shows that
the compact set ${\rm clos} \tn{\ve}{\cal A}$ is weakly invariant
for \rref{ext3}.
 Since  \rref{ext3} is GAC to ${\rm clos}
\tn{\ve}{\cal A}$ on $\R^k$, the CLSS Theorem (namely, Theorem
\ref{babyCLSS} above)
 provides an s-stabilizing feedback $K(x)$ for \rref{ext3}.
By Lemma \ref{invariance}, $\maf$ is strongly invariant for $f_2$.
 It follows that  the
$u$-part of $K(x)$ stabilizes~\rref{systema}.  This establishes
Theorem \ref{CLSS}.

\section{Illustration}
\label{illus} We next illustrate our stabilization approach using
the system
\begin{equation}
\label{spheredyn} \dot x=A_1(x)u_1+A_2(x)u_2\; \in\;  T_x(\maf),
\; \; x\in \maf, \; \; u=\langle u_1,u_2\rangle\in \R^2
\end{equation}
evolving on the sphere $\maf:=S^2=\{x\in \R^3: |x|=1\}$.  This
simple example  will illustrate how to construct  stabilizing
state feedbacks and
 Lyapunov functions on smooth manifolds. Even in this simple case,  we will see the
 necessity for using {\em discontinuous} stabilizers. Our example is a
modified version of the engineering examples in \cite{BMS95}. We
choose
 the
attractor ${\cal A}=\{\pm q\}$, where $q=\langle 0,0,1 \rangle $,
but similar constructions apply for any $q\in S^2$. The vector
fields $A_1$ and $A_2$ are chosen as follows. First define
$B_1(x)=q-(x\cdot q)x$ and $B_2(x)=x\times q$, which form an
orthogonal basis for the tangent spaces $T_x(\maf)={\rm
span}\{x\}^\perp$ on $S^2\setminus{\cal A}$ (in terms of the cross
product $\times$, the standard inner product $\cdot$,  and the
orthogonal complement $\perp$). Define the {\em geodesic distance}
${\cal G}$ on $S^2$ by
\[{\cal G}(x, x^\prime):={\rm arccos} (x\cdot x^\prime),\; \; {\rm
for}\;  x,x^\prime\in S^2.\] Set $r=\langle 0,1,0\rangle\in S^2$
and
\begin{equation}\label{vdef}V_q(x)=\min\{{\cal G}(x,\bar q): \bar q\in \{\pm q\}\}\;
  \; {\rm and}\; \;     V_r(x)=\max\{{\cal G}(x,\bar r): \bar r\in \{\pm
r\}\},\; \; \; x\in S^2.\end{equation}  Note the asymmetry between
$V_q$ and $V_r$.  Roughly speaking, we use a max in $V_r$ to
produce a component in our Lyapunov function that penalizes states
near $\pm r$ (see \rref{lya}). Let $M_1:S^2\to[0,1]$ be any smooth
function satisfying:
\begin{itemize}
\item $M^{-1}_1(0)=\left\{x\in S^2: \frac{x_1}{4}\le x_2\le
\frac{3\, x_1}{4}\;
 {\rm and}\;   V_q(x)\ge \frac{\pi}{4}\right\}$
\item$M^{-1}_1(1)\supseteq\left\{x\in S^2: x_2\ge \frac{7
x_1}{8}\; {\rm or} \; x_2\le \frac{x_1}{8} \;  {\rm or}\;
V_q(x)\le \frac{\pi}{8}\right\}$\end{itemize} and set
\begin{equation}
\label{vfs} A_1(x)=M_1(x)B_1(x),\; \;  A_2(x)=B_2(x).
\end{equation}
The factor $M_1$ in \rref{vfs} introduces a set of zeros in $A_1$,
consisting of a geodesic rectangle covering a part of the equator
of $S^2$ in the quadrant $Q_{++}:=\{x\in S^2: x_1>0, x_2>0\}$.  In
particular, the system \rref{spheredyn} is not completely
controllable. The fact that \rref{spheredyn} is GAC to ${\cal A}$
follows because any initial value can be moved to ${\cal A}$ along
the geodesic direction (i.e., ``north'' or ``south'' along a great
circle through $\pm q$) using the vector field $A_1$, possibly by
first using $A_2$ to move  the state ``west'' out of
$M^{-1}_1([0,1))$; see below for a precise definition of these
stabilizing trajectories. In fact, this global stabilization
is
done by the discontinuous feedback \rref{sfb} we construct
below.   On the other hand, a simple continuous dependence and
separation argument (e.g., the argument from the appendix in
\cite{S79})
 shows that the system has no Lipschitz
stabilizing state feedback $K(x)$.

The extension of \rref{spheredyn} from the Generalized CLSS
Theorem amounts to projecting onto the sphere, as follows. The
state space $\maf$ embeds into $\R^3$ by inclusion and
$\tube(\maf)=\{\langle x,kx\rangle: x\in \maf, k\in \R\}$. We can
choose $\omega(x)\equiv 1/4$ on $\maf$.  This  gives the
$\omega$-tube and annular tubular neighborhood
\[
\tube(\maf, \omega)=\{\langle x,kx\rangle\in \maf\times \R^3:
|k|<1/4\},\; \; \; \; \tn{\omega}{\maf}=\theta(\tube(\maf,
\omega))= \left\{x\in \R^3: 3/4<|x|<5/4\right\}.\] In terms of the
projection $\Pi_s(y):=y/|y|$ defined on $\R^3\setminus\{0\}$, our
corresponding  system $f_1$ on $\tn{\omega}{\maf}$ is (see
\rref{ext2})
\[
\null\hskip-0.8cm
f_1(y,u,v)=M_1(\Pi_s(y))\left(q-\left\{\Pi_s(y)\cdot
q\right\}\Pi_s(y)\right) u_1+\left\{\Pi_s(y)\times q\right\}u_2+
\left(y-\Pi_s(y)\right)v,\; \; \langle u,v\rangle\in \R^2\times
\R,
\]
which we (discontinuously) extend to $\R^3$ by setting $f_1\equiv
0$ for states outside $\tn{\omega}{\maf}$.
 We next choose
\[\maf^\sharp=\{z\in \R^3: 7/8\le |z|\le 9/8\},\; \; \;  C_\omega=\R^3\setminus \{z\in
\R^3: 13/16\le |z|\le 19/16\}.\]
 Our corresponding system
$f_2$ on $\R^3$
 can then be defined by
taking $\phi(z)=\Gamma(|z|^2)$ for any smooth  function
$\Gamma:[0,\infty)\to [0,1]$ that satisfies (i) $\Gamma\equiv 1$
on $[(7/8)^2, (9/8)^2]$ and (ii) $\Gamma\equiv 0$ outside
$((13/16)^2, (19/16)^2)$.  Since $f_2$ is GAC to $\tn{\omega/{\rm
2}}{{\cal A}}$, there exists a sample stabilizing feedback for
$f_2$ whose restriction $K(x)$ to $\maf$ stabilizes
\rref{spheredyn} to ${\cal A}$.  This is the content of our
Generalized CLSS Theorem.

The stabilizing feedback $K(x)$  and a corresponding
control-Lyapunov function (CLF) can be explicitly constructed by
the following variant of the argument from \cite[Section
2]{BMS95}. Set ${\cal A}^\sharp=\{\pm q, \pm r\}\subseteq S^2.$
Define
\[Y^{\bar p}_x:= \Pi_s(\bar p-\{\bar p\cdot x\}x),\; \; \bar p\in
{\cal A}^\sharp,\; x\ne \pm \bar p;\] this
 gives the geodesic direction from $x$
to $\bar p$.  Note that \[A_1(x)\cdot Y^q_x\equiv \sqrt{1-x^2_3}\;
\;  {\rm and} \; \; A_2(x)\cdot Y^q_x \equiv 0.\]  Also, $Y^{-\bar
p}_x\equiv -Y^{\bar p}_x$ for all $\bar p\in {\cal A}^\sharp$.
 A straightforward calculation (see \cite[Lemma 1]{BMS95}) shows
 that along any (open loop) trajectory of
\rref{spheredyn}  that does not pass through ${\cal A}^\sharp$, we
get
\begin{equation}
\label{gauss} \frac{d}{dt} {\cal G}(x,\bar p)=- \dot x\cdot
Y^{\bar p}_x,\; \; \bar p\in {\cal A}^\sharp,\; x\ne \pm \bar p.
\end{equation}
We  show that \rref{spheredyn} can be globally stabilized to
${\cal A}$ by the (necessarily discontinuous) state feedback
\begin{equation}\label{sfb}
K(x)= \left\{
\begin{array}{rl}
\langle 1,0\rangle, & {\rm if}\; \;  x_3\ge 0\; \; {\rm and}\; \;   M_1(x)=1\\
-\langle 1,0\rangle, & {\rm if}\; \; x_3<0\; \; {\rm and}\; \;   M_1(x)=1\\
\langle 0,1\rangle, & {\rm if}\; \; M_1(x)<1
\end{array}
\right.
\end{equation}
when the closed-loop trajectories are defined in the usual
non-sampling sense.  An easy argument will then show that
\rref{sfb} also sample stabilizes \rref{spheredyn}.
 Before presenting our argument, we
interpret \rref{sfb} in terms of the corresponding closed loop
(non-sampling) trajectories. For values where $M_1(x)=1$, the
feedback $K$ drives the state to ${\cal A}$ geodesically along a
great circle through $\pm q$. On the other hand, any state where
$M_1(x)<1$ is driven towards $- r$ until the state reaches
$M^{-1}_1(1)$ and then geodesically to ${\cal A}$.

We first analyze the usual non-sampling trajectories of the closed
loop system for $K(x)$ which we refer to simply as ``closed loop
trajectories'' in the sequel.
 The fact that $K(x)$ stabilizes the closed loop trajectories to ${\cal A}$
can be verified using the following Lyapunov function
construction. In terms of  $V_q$ and $V_r$ in \rref{vdef}, set
\begin{equation}\label{lya} V(x)=V_q(x)[1+V_r(x)],\; \; x\in S^2.\end{equation}
 Then $V$ is continuous and nonnnegative, and
$V$ is null only on ${\cal A}$.  We will show that $V$ is an {\em
integral} Lyapunov function for \rref{spheredyn} in the sense of
\cite{AS99}; this will  imply that $V$ is also a CLF in the usual
Dini derivative sense used for example in
\cite{CLSS97}.\footnote{A {\em control-Lyapunov integral function}
for \rref{spheredyn} and ${\cal A}$ is defined to be any
continuous function $V:\maf\to [0,\infty)$ for which
$V^{-1}(0)={\cal A}$ and for which there exist a constant $N>0$
and $\alpha_3\in {\cal K}$ satisfying:  For each $\xi\in \maf$,
there exists $u\in {\cal U}_N$ such that $x(t):=x(t,\xi,u)$ is
well defined and satisfies (see \cite{AS99})
\begin{equation}\label{starr}
V(x(t))-V(\xi)\le -\int_0^t \alpha_3(|x(s)|_{\cal A})ds \; \;
\forall t\ge 0.\end{equation} We will verify the decay condition
\rref{starr} using closed loop trajectories and corresponding
feedback controls $u(t)=K(x(t))$ for \rref{spheredyn}. The
inequality \rref{starr} then gives the usual Dini derivative
Lyapunov decay condition for $V$ (e.g., from \cite{CLSS97, KT04})
once we divide through by $t$ and pass to the liminf.  This last
step uses the fact that $t\mapsto \dot x(t)=f(x(t),K(x(t)))$ is
(right) continuous at $t=0$ for each closed loop trajectory $x(t)$
of \rref{spheredyn}.}
 Given a closed loop trajectory $x(t)$, we also let $\dot V(x)$ denote the
derivative of $t\mapsto V(x(t))$ when it is defined. Along any
trajectory $x(t)$ of the closed loop system that remains in
$M^{-1}_1(1)\setminus {\cal A}^\sharp$ and that satisfies $x_2>0$
and $x_3>0$ everywhere, we get $V(x)={\cal G}(x,q)[1+{\cal
G}(x,-r)]$ and therefore \rref{gauss} gives
\begin{eqnarray}
\dot V(x)& =&  -\dot x\cdot Y^q_x\, [1+{\cal G}(x,-r)]-\dot x\cdot
Y^{-r}_x\, {\cal G}(x,q)\nonumber\\  & =& -A_1(x) \cdot Y^q_x\,
[1+{\cal
G}(x,-r)]+A_1(x)\cdot\Pi_s(r-\{x\cdot r\}x){\cal G}(x,q)\nonumber\\
&=&-A_1(x)\cdot Y^q_x[1+{\cal G}(x,-r)] -\frac{x_2\,
x_3}{|r-(x\cdot r)x|}\, {\cal G}(x,q)\; \; \le \; \;
-\sqrt{1-(x\cdot q)^2}\nonumber,
\end{eqnarray}
which is only zero when $x\in {\cal A}$. Similar arguments show
that
\begin{equation}
\label{decay1} \dot V(x)\le -\sqrt{1-(x\cdot q)^2},\; \; {\rm
if}\; \;  x_2\ne 0,\; \; x_3\ne 0,\; \; M_1(x)=1,\; \; {\rm and}\;
\; x\not\in {\cal A},
\end{equation}
and  $\dot V(x)=-(1+\frac{\pi}{2})\sqrt{1-x^2_3}$ along
trajectories in $\{x_2=0, x_3\ne 0\}$.  Notice that $\dot V$ is
continuous along closed loop trajectories in $M^{-1}_1(1)$
starting outside $\{x_3=0\}$, and that the closed loop
trajectories  starting in $M^{-1}_1(1)$ with $x_3(0)=0$ also
satisfy $x_3(t)\ne 0$ for all $t>0$. This gives
$V(x(t))-V(x(0))\le -\int_0^t \sqrt{1-x^2_3(s)}ds$ along each
closed loop trajectory starting in $M^{-1}_1(1)$.

On the other hand, along closed loop trajectories in
$M^{-1}_1([0,1))$, we know that $V_r(x)={\cal G}(x,-r)$, so
\begin{eqnarray}\label{decay2}
\dot V(x)& =& -\dot x\cdot Y^{\pm q}_x\, [1+{\cal G}(x,-r)]-\dot
x\cdot Y^{-r}_x\, V_q(x)\; =\;    -A_2(x)\cdot Y^{-r}_x\, V_q(x)\;
\;  \; ({\rm since\ } A_2(x)\cdot Y^q_x\equiv 0)
\nonumber\\
&=& (x\times q)\cdot \Pi_s(r-\{x\cdot r\}x)V_q(x)\; \; =\; \;
-\frac{x_1\, V_q(x)}{\sqrt{x^2_1+x^2_3}}\; \; =:\; \;
-\mu(x)\end{eqnarray} when $x_3\ne 0$.  Notice that $-\mu$
 is bounded above by a negative
constant in $M^{-1}_1([0,1))$. Also, $\dot V(x)\equiv -\pi/2$
along closed loop trajectories in $M^{-1}_1([0,1))$ along $x_3=0$.
Therefore,  reasoning exactly as before gives \[V(x(t))-V(x(0))\le
-\int_0^t \mu(x(s))ds\] along all closed loop trajectories $x(t)$
remaining in $M^{-1}_1([0,1))$.  Since $M^{-1}_1(1)$ is forward
invariant for the closed loop trajectories, it  follows that
 the discontinuous feedback $K(x)$
stabilizes  the closed loop trajectories of \rref{spheredyn} to
${\cal A}$, and that $V$ satisfies the requirements for being a
control-Lyapunov (integral) function for \rref{spheredyn} and also
 a CLF for \rref{spheredyn} in the usual Dini derivative sense of
\cite{CLSS97}.  That $K(x)$ also {\em sample} stabilizes
\rref{spheredyn} now follows because (a) the sampling and
(non-sampling) closed loop trajectories agree for initial points
in $M^{-1}_1(1)$ and (b) the equality $\dot V(x)=-\mu(x)$ holds
throughout the quadrant $Q_{++}$ if we use the control $u\equiv
\langle 0,1\rangle$ at all points in  $Q_{++}$.  In fact, (a)
implies that $K$ sample stabilizes \rref{spheredyn} for initial
values in $M^{-1}_1(1)$ for all partitions $\pi$.  Also, (b)
implies that $K$ sample stabilizes the dynamic for initial values
in $M^{-1}_1([0,1))$ when  $\overline{\mathbf d}(\pi)$ is
sufficiently small for the sample control value to switch to $\pm
\langle 1,0\rangle$ in $M^{-1}_1(1)$ but before the first time the
sample trajectory exits $Q_{++}$.

\section{Further Extensions}
\label{sec5}

We next use our results   to establish the input-to-state
stabilizability (ISSability) of control affine systems
\begin{equation}
\label{affine} \dot x=h(x)+G(x)u,\; \; x\in {\cal X},\; \; u\in
{\mathbf U}=\R^m
\end{equation}
evolving on smooth Riemannian manifolds ${\cal X}$ relative to
actuator errors  (but see Remark \ref{rk1} below for an extension
to fully nonlinear systems).  We assume \rref{affine} is GAC to a
weakly invariant compact nonempty set ${\cal A}\subseteq\maf$. In
this context, $G(x)u= g_1(x)u_1+\ldots +g_m(x)u_m$ for locally
Lipschitz vector fields $g_i:{\cal X}\to\R$.
 The stabilizers we construct in this
section have the additional desirable feature that they are robust
to small observation noise in the controllers.  For  {\em
continuous} feedback stabilizers, small observation noise in the
controllers can be tolerated.  However, since  our stabilizing
feedback may need to be  {\em dis\/}continuous (see Section
\ref{sec1}), such noise terms can have a substantial effect on the
dynamics.  Therefore,  the magnitude of the noise needs to be
constrained in terms of the sampling frequency (see \cite{MRS04,
MS04, S99} and Definition \ref{sampleISSstab} below).

  To make
our ISSability notion precise, we first introduce a Riemannian
metric on ${\cal X}$ to quantify  observation noise and we let
${\cal B}_r(y)$ denote the corresponding closed ball in ${\cal X}$
centered at $y\in {\cal X}$ of radius $r$.  As before, a {\em
feedback} for (\ref{affine}) is  defined to be any locally bounded
function $k:\maf\to {\mathbf U}$.  We introduce the set of
functions ${\cal O}=\{e:[0,\infty)\to [0,\infty)\}$, which
represent the observation errors in our controller; and for each
$e\in {\cal O}$, we set $\sup(e)=\sup\{e(t): t\ge 0\}$.  We use
the set of functions ${\cal O}(\eta):=\{e\in {\cal O}: \sup(e)\le
\eta\}$ for each $\eta>0$. We let ${\rm Par}$ denote the set of
all partitions and \[{\rm Par}(\delta):=\left\{\pi\in {\rm Par}:
\overline{\mathbf d}(\pi)<\delta\right\}\] for each $\delta>0$.
Our ISSability goal of this section is to find a feedback $k$ so
that
\begin{eqnarray}
\dot x(t)&=&h(x(t))+G(x(t))[k(\eta(t))+\uc(t)],\; \; \eta(t)\in
{\cal B}_{e(t)}(x(t))
\label{l2}
\end{eqnarray}
is input-to-state stable (ISS) for sampling solutions relative to
actuator errors $\uc$ for small  observation errors $e$.  The
relevant definitions are as follows:

\begin{definition}
\label{sampleISSsol} Let $k$ be a feedback for (\ref{affine}),
$e\in {\cal O}$, $\uc\in {\cal U}$, $\xi\in {\cal X}$, and
$\pi=\{t_i\}_{i\ge 0}$ be any partition of $\R_{\ge 0}$.  A {\em
$\pi$-solution} for (\ref{l2}), the initial state $\xi$, the
observation error $e\in {\cal O}$, and ${\mathbf u}$ is defined to
be any continuous function $x(\cdot)$ obtained by recursively
choosing any $\eta(t_i)\in {\cal B}_{e(t_i)}(x(t_i))$ and then
solving
\[
\dot x(t)=h(x(t))+G(x(t))[k(\eta(t_i)) +\uc(t)]
\]
from the initial  time $t=t_i$ up to time
\begin{equation}\label{sii}s_i=\max\{t_i,\sup\{s\in [t_i,t_{i+1}]: x(\cdot) {\rm\ is\
defined\ on\ } [t_i,s)\}\},\end{equation} where
$x(0)=\xi$.\footnote{As before, the continuity requirement for
$x(\cdot)$ stipulates that the final value on the previous
subinterval is used as the initial value at the next subinterval.
Also, the $t_i$ argument of the max \rref{sii} allows the
possibility that $x(t)$ is not defined at all on $[t_i,t_{i+1}]$
(see Section \ref{sec2}). }
 The domain of
$x(\cdot)$ is $[0,t_{\rm max})$, where $t_{\rm max}=\inf\{s_i:
s_i<t_{i+1}\}$. When $t_{\rm max}=+\infty$, we call $x(t)$ {\em
well defined}. \mybox
\end{definition}

\begin{definition}
\label{sampleISSstab} Let $k$ be a feedback for (\ref{affine}). We
say $k$ renders (\ref{affine}) {\em sample-input-to-state stable
(s-ISS) to ${\cal A}$}
 provided for each ${\cal R}_o\in {\cal PN}^{\cal A}$ and each
$N>0$, there exists ${\cal R}_1={\cal R}_1(N) \in {\cal PN}^{{\cal
R}_0}$ such that:
\begin{enumerate}\item
For each ${\cal R}_2,{\cal R}_3\in {\cal PN}^{{\cal R}_1}$ with
${\cal R}_2\subseteq {\cal R}_3$, there exist ${\cal M}={\cal
M}({\cal R}_3)\subseteq {\cal X}$ and positive numbers $\delta$,
$T$, and $\kappa$ (depending  on ${\cal R}_2$ and ${\cal R}_3$)
such that if $\pi\in {\rm Par}(\delta)$, $\xi\in {\cal R}_3$,
$\uc\in {\mathbf {\cal U}}_N$, and $e\in {\cal O}(\kappa
\underline{\mathbf d}(\pi))$, then the corresponding
$\pi$-solutions $x(t)$   for (\ref{l2}) starting at $\xi$ are  all
well defined and satisfy (a) $x(t)\in {\cal R}_2$ for all $t\ge T$
and (b) $x(t)\in {\cal M}$ for all $t\ge 0$. \item For each ${\cal
E}\in {\cal PN}^{{\cal R}_1}$, there exists $D\in{\cal PN}^{{\cal
R}_1}$ such that if the set ${\cal R}_3$ in 1. is a subset of $D$,
then the set ${\cal M}$ in 1. can be chosen to be a subset of
${\cal E}$.
\end{enumerate}
and for each ${\cal R}_0\in {\cal PN}^{\cal A}$, there exists
$N=N({\cal R}_0)>0$ such that 1.-2. hold with the choices
$N=N({\cal R}_0)$ and ${\cal R}_1={\cal R}_0$. In this case, we
also say (\ref{l2}) is {\em ISS for sampling solutions} and that
\rref{affine} is {\em ISSable}.\mybox
\end{definition}

The preceding definition requires that the sampling  be done
quickly enough so that $\pi\in {\rm Par}(\delta)$, but not so
quickly that ${\rm sup}(e)>\kappa \underline{\mathbf d}(\pi)$.
When  $e\equiv 0$, the condition on $\underline{\mathbf d}(\pi)$
in  Definition \ref{sampleISSstab} is not needed.  For ${\cal
X}=\R^n$, one can easily check that if \rref{affine} is sampling
ISS in the sense defined in \cite{MRS04, MS04} using some feedback
$k$, then it is also ISSable in the sense of Definition
\ref{sampleISSstab} with the same feedback $k$. For any compact
nonempty weakly invariant set ${\cal A}$ for \rref{affine}, we
then have:

\bt{affinethm} If (\ref{affine}) is GAC to ${\cal A}$, then there
exists a feedback $k(x)$ rendering \rref{affine} s-ISS to ${\cal
A}$. \et

\begin{proof}
We indicate the changes needed in the proof of Theorem \ref{CLSS}.
As before, we first extend the dynamics \[f(x,u)=h(x)+G(x)u\] to a
dynamics \rref{ext3} defined on all of $\R^k$ that is GAC to ${\rm
clos}\tn{\ve}{\cal A}$.  By \cite[Theorem 3.2]{KT04}, this
extended dynamics admits a locally Lipschitz control-Lyapunov
function (CLF) $V$; see \cite{S99} for  background on CLFs. Using
the  argument from \cite[Section 5]{R00}, we can transform $V$
into a (locally) semiconcave CLF for \rref{ext3}  on
$\R^k\setminus{\rm clos}\tn{\ve}{\cal A}$. In \cite{MRS04}, it was
shown that control affine systems that are GAC to $\{ 0\}$ on
$\R^k$ admit (possibly discontinuous) feedbacks for which the
corresponding closed loop systems are sampling ISS to $\{0\}$.
Since  \rref{ext3} is again control affine, a slight variant of
the argument from \cite[Section 3]{MRS04}
provides a feedback $K(x)$ rendering \rref{ext3} s-ISS to ${\rm
clos}\tn{\ve}\A$.  Applying  Lemma \ref{invariance} as before, we
conclude that the $u$-part of $K(x)$  renders \rref{affine} s-ISS
to ${\cal A}$.
\end{proof}

\begin{remark}
\label{rk1} The preceding theorem can be extended  to cover fully
nonlinear systems \rref{systema} on $\maf\times \R^m$ if we
reinterpret s-ISS in the following more general sense:  A feedback
$k$ renders \rref{systema} {\em s-ISS to ${\cal A}$ in the weak
sense} provided there exists a smooth everywhere invertible matrix
valued function $G:{\cal X}\to \R^{m\times m}$ such that \be{weak}
\dot x=f(x,k(x)+G(x)u) \ee is s-ISS to ${\cal A}$.  The s-ISS
property for \rref{weak} is defined by taking $e\equiv 0$ in
Definition \ref{sampleISSstab}, and the $\pi$-solutions of
\rref{weak} are defined by recursively solving \[\dot
x(t)=f(x(t),k(x(t_i))+G(x(t))\uc(t))\] on successive intervals
$[t_i, t_{i+1}]$ of the partition $\pi=\{t_i\}_{i\ge 0}$ and
proceeding as in  Definition \ref{sampleISSsol} with $e\equiv 0$
(see \cite{MRS04} for details).  In particular, the sampling is
only done in the (possibly discontinuous) controller $k(x)$. We
can then prove the following for any smooth manifold ${\cal X}$
and ${\mathbf U}=\R^m$: {\em If \rref{systema} is GAC to a
compact, nonempty, weakly invariant set ${\cal A}$, then there
exists a feedback $k(x)$ rendering \rref{systema} s-ISS to ${\cal
A}$ in the weak sense.} The proof combines the arguments from
\cite[Section 5]{MRS04} with our proof of Theorem \ref{affinethm}
and is left to the reader.\end{remark}

\begin{remark}
\label{rk3} As we noted in the introduction, the GAC system
\rref{systema} will not in general admit a continuous stabilizing
state feedback. However,  by \cite{C92}, the system \rref{systema}
is stabilizable by a continuous {\em time varying} feedback
$u=k(t,x)$ if it is completely controllable and drift-free (the
latter condition being the requirement that $f(x,0)\equiv 0$). In
engineering applications, feedback laws are  usually implemented
via sampling.  This motivated our construction of discontinuous
state stabilizers $u=k(x)$ which we implemented using CLSS
solutions. Yet another approach to stabilizing \rref{systema} is
to look for a {\em dynamic} stabilizer. This means finding a
locally Lipschitz regulator dynamic $\dot z=A(z,x)$ and a locally
Lipschitz function $k(z,x)$ such that the interconnected system
\[\dot x=f(x,k(z,x)),\; \;  \dot z=A(z,x)\] is globally asymptotically
stable.  See \cite{S98} for an extensive discussion of dynamic
stabilizers for {\em linear} systems.

On the other hand, it turns out that a dynamic feedback for
\rref{systema} may fail to exist, even if the system is completely
controllable. An example from \cite{S79} where this occurs is
\begin{equation}
\label{hjs}
 \dot x=f(x,u)= \left[
\begin{array}{l}
(4-x^2_2)u^2_2\\
e^{-x_1}+x_2-2e^{-x_1}\sin^2(u_1)
\end{array}
\right],\; \; x\in \R^2,\; u\in \R^2.
\end{equation}
The fact that \rref{hjs} is completely controllable (and therefore
GAC to ${\cal A}=\{0\}$) was shown in the appendix of \cite{S79},
where it is also shown that it is impossible to choose paths
converging to the origin in such a manner
  that this selection is continuous in the initial states.
  Since the flow map of any dynamic
stabilizer would give a continuous choice of paths converging to
the origin, no dynamic stabilizer for the system can exist, even
if we drop the requirement that the state of the regulator
converges to zero. In particular, we see that \rref{hjs} cannot
admit a continuous time varying feedback $u=k(t,x)$.  This does
not contradict the existence results \cite{C92} for time varying
feedbacks since in this case, the system has drift.
\end{remark}
\begin{remark}
The feedback construction  \cite{MRS04} used to prove Theorem
\ref{affinethm} proceeds by first finding a semiconcave
control-Lyapunov function (CLF) for the system and then adapting
the feedback design from \cite{S89} to allow nonsmooth CLFs,
observation noise,  and discontinuous feedback. Semiconcave CLFs
are known to exist for all (locally Lipschitz) GAC systems on
Euclidean space and all compact nonempty weakly invariant
attractors ${\cal A}$,  by arguments from \cite{R00}. The
semiconcavity property is intermediate between $C^1$ and local
Lipschitzness. On the other hand, GAC systems will not in general
admit {\em smooth} CLFs since their existence would imply the
existence of continuous stabilizers $k(x)$, which we know not to
be the case in general (see \cite{B83, S99}).

For a very different approach to ISS on manifolds (based on
density functions) that gives rise to a sufficient condition for
ISS-like behavior from almost all initial values, see \cite{A04}.
The main ISS-like condition in \cite{A04} states: For a given
Riemannian manifold ${\cal X}$ and a compact weakly invariant set
${\cal A}\subseteq {\cal X}$ for \rref{systema}, we say that
\rref{systema} is {\em weakly almost ISS to ${\cal A}$} provided
(i) ${\cal A}$ is locally asymptotically stable for the system and
(ii) there exists $\gamma\in {\cal K}$ such that
\begin{equation}
\label{DAcond} \forall u\in {\cal U},\; \; \exists{\cal Z}_u\in
{\rm Null}({\cal X})\; \; {\rm s.t.}\; \; \forall \xi\in {\cal
X}\setminus {\cal Z}_u, \; \; \liminf_{t\to
+\infty}|x(t,\xi,u)|_{{\cal A}}\le \gamma({||u||})
\end{equation}
where ${\rm Null}({\cal X})$ is the set of subsets of $\maf$ of
measure zero and $|\cdot|_{{\cal A}}$ denotes the distance to
${\cal A}$.
 This condition differs from our ISS requirement mainly in its
allowance of a null set of states that are not necessarily
stabilized and in its use of Carath\'eodory solutions.
An alternative and more intrinsic approach to feedback
stabilization on manifolds would involve generalizing the concepts
of set valued differentials and semiconcave CLFs to manifolds and
providing direct feedback constructions without first embedding
into $\R^k$. We provided a first result in this direction in
Section \ref{illus} above. We leave the development of this more
intrinsic approach for another paper.
\end{remark}
\newpage
\section*{Acknowledgments}
  M. Malisoff was supported
 by  Louisiana Board of Regents  Contract
 LEQSF(2003-06)-RD-A-12.  Part of the work of M. Krichman  was carried out while this author
was a Research Assistant at  Rutgers University. Krichman thanks
Felipe M. Pait for helpful comments.  E. Sontag was supported by
NSF Grant CCR-0206789.

\newpage


\end{document}